\documentclass{article}

\usepackage[dblblindworkshop, final]{neurips_2026}

\usepackage[utf8]{inputenc} 
\usepackage[T1]{fontenc}    
\usepackage{hyperref}       
\usepackage{url}            
\usepackage{booktabs}       
\usepackage{amsfonts}       
\usepackage{nicefrac}       
\usepackage{microtype}      
\usepackage{xcolor}    

\newtheorem{definition}{Definition}[section]

\newtheorem{theorem}{Theorem}[section]
\newtheorem{lemma}{Lemma}[section]

\title{Proof Perplexity: A Probabilistic Approach to Measuring non-triviality in Machine Generated Proofs}

\author{%
Madhuparna Das\thanks{The author has since moved to a different position and is submitting this work independently.}\\
Department of Computer Science and Technology\\
University of Cambridge\\
Cambridge, UK\\
\texttt{md2147@cam.ac.uk} \\
\AND
Soumya Banerjee\\
Department of Computer Science and Technology\\
University of Cambridge\\
Cambridge, UK\\
  \texttt{sb2333@cam.ac.uk} 
}

\begin{document}

\maketitle

\begin{abstract}
In this article, we propose a benchmark to quantify the non-triviality and novelty of machine-generated mathematical proofs compared to existing literature. Our core methodology uses a probabilistic measure to evaluate the non-triviality of the underlying mathematical arguments in comparison to the current \texttt{MATHLIB} library. Finally, we apply this benchmark to the recent proof generated by Claude on the complex zeros of the Riemann zeta function to test for non-trivial steps.
\end{abstract}

\section{Introduction}

The recent development of machine-assisted provers for mathematics has been growing significantly, especially with the emergence of generative AI. While it is difficult for LLMs to verify whether their own generated proof is mathematically correct or not, interactive theorem provers such as Lean and Isabelle satisfy this step quite well. Although these systems have manually curated and maintained libraries, they nonetheless work as an automated verification layer for mathematical proofs. As LLMs are now solving mathematical research problems that can be easily verified by formalising them in systems like Lean. One important question still remains:

\noindent
{\it How much of machine-generated proofs are truly novel discoveries rather than the use of known results from the literature?}

\noindent
Since mathematics is a language, it is sometimes subjective for a proof step to be labelled as non-trivial. For example, one might argue that a particular proof step is a mere recombination of known results from the literature and that the machine has not discovered something novel. However, for the same proof step, one might say that such recombination was not immediately observed by human mathematicians and opened new frontiers for further research, and therefore is non-trivial. Thus, it is a challenge to properly measure whether a machine-generated proof is truly novel and non-trivial or a mere reconstruction of known results. In this article, we establish a methodology to probabilistically measure the non-triviality of a machine-generated proof and use Claude's proof of the distribution of the zeros of the Riemann zeta function as an example for our experiment. We use a probabilistic method because it allows us to see if a proof step is a genuine statistical outlier compared to existing mathematical literature within formal libraries, which gives us a much more accurate measure of the novelty of the discovery than simply assuming a proof step is novel because it is not symbolically similar to the known results.

{\it Related Results.} N. Somani's Prior Proof~[1] measures proof-route novelty against historical \texttt{MATHLIB} snapshots, capturing how a proof’s novelty changes as similar dependency patterns become more common over time. Similar work on studying the novelty and non-triviality of machine-generated proofs was done by [2, 3].

\subsection{Notations and Initial Setup}\label{1.1}

To measure how much of a machine-generated proof is genuinely novel, rather than simply retrieved or reformulated from results already available in a formal library, we need to examine how automated provers interact with existing mathematical libraries. Automated systems such as AlphaProof~[8] and Aristotle~[9] use search algorithms that naturally represent proof exploration as a tree. Because this search space can grow rapidly, we need to define a clear boundary around the existing library to determine which parts of a generated proof constitute genuinely new results. It is worth noting that our proposed methodology can be used for machine-generated proofs that are formalised in the \texttt{MATHLIB} library. We use the following notations throughout this article.

Let $\mathcal{L}$ be the entire \texttt{MATHLIB} library without our fixed target theorem $T$ and its proof being formalised. Assume that the target theorem has $n\in\mathbb{N}$ intermediate declarations or proof steps, and let these proof steps be denoted by $s_i$, for $1\leq i\leq n$. These proof steps can be lemmas, propositions, or theorems, and can have their own proofs. For a fixed theorem $T$, we consider its dependency closure, i.e., the set of Lean declarations on which the proof steps of $T$, as checked by the Lean kernel, actually depend. We then arrange these proof steps in a topological order, denoted by $\pi=(s_1,\ldots,s_n)$, such that each $s_i$ depends only on declarations in $\mathcal{L}\cup\{s_1,\ldots,s_{i-1}\}$. We define the internal complexity $w(s_i)$, which denotes the internal declarations of the proof step $s_i$. 

\begin{definition}[Search Budget]\label{searchb}
Let $\mathcal{T}(s_i)$ denote the proof search tree generated by an automated prover with the target step $s_i$ as its root node. The search budget, denoted by $B \in \mathbb{N}$, is defined as the maximum allowable depth of \;$\mathcal{T}(s_i)$. For any valid search path $p \in \mathcal{T}(s_i)$, the path length $|p|$ is bounded strictly by $B$, such that $|p|\leq B$. Search branches exceeding depth $B$ are terminated\footnote{Note that one may choose the time parameter or tokens as their search budget rather than the size of the tree path. }.    
\end{definition}

Let $w(s_i, B)$ denote the total number of tactic applications in the machine-generated proof for the specific target step $s_i$ under the search budget $B$. The prover is restricted from using the target step itself or any subsequent dependencies, relying strictly on $\mathcal{L}\cup\{s_1, \dots, s_{i-1}\}$. If the prover fails to find a valid proof within the budget $B$, we set $w(s_i, B):= w(s_i)$. Finally, we define the candidate pool, denoted as $\mathrm{Cand}_B$, which is a curated shortlist of declarations taken from the main library $\mathcal{L}$. Automated provers build this pool automatically using a smart filtering method to pick out the most relevant results for proving the target theorem $T$. 

For the experimental purpose of our methodology, we use Claude's new proof~[4] that $\frac{2}{3}$ of the zeros of the Riemann zeta function $\zeta(s)$ lie on the $\Re(s)=\frac{1}{2}$, and refer to it as the zeta23 proof throughout this article. Note that the proof has already been formalised in \texttt{MATHLIB} by Anthropic~[15].

\section{The Methodology}

In information theory and machine learning, \emph{perplexity} evaluates how well a probability model predicts a sample~[5]. For natural language processing, it measures a model's uncertainty or \emph{surprise} when predicting the next word in a sequence. A lower perplexity indicates that the model is highly confident and the generated text is mathematically predictable, while a higher perplexity means the model is confused and the sequence is unexpected.

Our methodology is motivated by this simple notion that effects the decision making process of a machine. We ask ourselves, if a step of a mathematical proof of the target theorem $T$ is a simple recombination or a direct extraction of the known results, then can it simply be regenerated or \emph{recovered} by the same machine with low surprise? If the observation is non-trivial, then the recovering of that proof step will not be so simple, assuming that the machine does not have the proof of target theorem $T$ in its memory. Motivated by this idea, we introduce a similar notion of perplexity to probabilistically measure the surprise of a proof step generated by a machine. However, measuring the surprise of a machine-generated proof step is not straightforward, as we first must measure the recoverability of a proof step.


\subsection{A Probabilistic Measure for Recoverability}

The real complexity in recovering a proof of a theorem is not recovering the step $s_i$ itself, rather generating the proof of the proof step $s_i$ for the target theorem $T$ within the search budget $B$ (see Definition~\ref{searchb}) from the library $\mathcal{L}$ because for the theorem $T$, $s_i$ lies in the first branch of the proof tree $\mathcal{T}(s_i)$, wherein the internal complexity $w(s_i)$ lies in the later branches, which makes it harder to be immediately recovered. To quantify this notion, we first introduce the concept of \emph{ablation recoverability}. By doing this, we evaluate the complexity of the newly generated (or recovered) proof, denoted as $w(s_i, B)$, against the internal complexity of the original step, $w(s_i)$.

\begin{definition}[Ablation recoverability]\label{ablationr}
For a given proof step $s_i$ and search budget $B$, the ablation recoverability score is defined as
\[
\alpha_{\mathrm{abl}}(s_i,B)
:=\max\left(1-\frac{w(s_i,B)}{w(s_i)},\ \alpha_{\min}\right),
\]
where the absolute minimum bound $\alpha_{\min}$ is defined as $\alpha_{\min}:=\frac{1}{1+|\mathrm{Cand}_B|}$.
\end{definition}

If $w(s_i, B)$ is small relative to $w(s_i)$, it implies the step is a standard, easily recoverable deduction. Conversely, if the prover exhausts the budget $B$ or requires a significantly deeper path, it indicates that the original step $s_i$ likely contained a non-trivial insight.  

\subsubsection{Mahalanobis Recoverability}\label{2.1.1}

The ablation recoverability captures the computational effort required to find a proof path for the proof tree; it is not sufficient on its own. A computationally difficult step such that the search budget is exhausted does not necessarily ensure that the underlying mathematical argument of that step is non-trivial or novel compared to the known literature. To address this problem, we treat the proof itself as a mathematical object and study its statistical distribution within the library $\mathcal{L}$. We define Mahalanobis recoverability~[12], which introduces a notion of correlational proxy that measures whether a machine-generated step $s_i$ belongs to the same statistical neighbourhood in the library $\mathcal{L}$ of the argument that came before it.

Let $\phi(s_i)$ be a structural feature vector for step $s_i$, a set of dependencies or premise usage vector over cited constants, or an embedding of the complexity $w(s_i)$ of the proof steps. Partition $\mathcal{L}$ into clusters $c=1,\dots,C$ by module or topic area, and fit a Ledoit-Wolf shrinkage~[6] $(\mu_c,\Sigma_c)$ from $\{\phi(d):d\in \mathcal{L},\ d\in c\}$ for each cluster. Define
\[
D_M(s_i\mid c)= \sqrt{(\phi(s_i)-\mu_c)^\top \Sigma_c^{-1}(\phi(s_i)-\mu_c)}.
\]
Under a multivariate $t$ null with $\nu_c$ degrees of freedom, fitted to each cluster $c$ alongside $(\mu_c, \Sigma_c)$ (see Appendix~\ref{appendixformah}), the rescaled quadratic form satisfies 
\[
\frac{\nu_c}{q(\nu_c-2)}\,D_M(s_i\mid c)^2\sim F_{q,\,\nu_c},\quad \mathrm{where}\; q=\dim\phi,
\]
giving a tail probability
\[
p_c(s_i) = 1 - F_{q,\,\nu_c}\left(\frac{\nu_c}{q(\nu_c - 2)}\,D_M(s_i\mid c)^2\right). 
\]
As $\nu_c \to \infty$, this recovers the Gaussian calibration $D_M(s_i\mid c)^2 \sim \chi^2_q$. The question is whether $s_i$ is structurally typical of at least one cluster, which is an intersection-union test. We take $\max_c p_c(s_i)$, if $s_i$ genuinely belongs to some cluster $c_0$, then $p_{c_0}(s_i)$ is stochastically at least $\texttt{Uniform}[0,1]$ under $H_{0,c_0}$, so $P_{H_0}(\max_c p_c < \delta) \leq P_{H_{0,c_0}}(p_{c_0}<\delta)\leq\delta$, controlling the type I error at level $\delta$ without multiplicity correction.

Using this maximum $p$-value to quantify how structurally typical a step is across the existing clusters, we formally define Mahalanobis recoverability.

\begin{definition}[Mahalanobis Recoverability]\label{mahr}
Let $\tilde{p}(s_i)$ denote the maximum $p$-value across all $C$ clusters
\[
\tilde{p}(s_i) 
:= \max_{1\leq c\leq C} p_c(s_i).
\]
The distributional recoverability of step $s_i$ is bounded below by $\alpha_{\min}$ and defined as
\[
\alpha_{\mathrm{mah}}(s_i) 
:= \max \big\{ \tilde{p}(s_i), \; \alpha_{\min} \big\}.  
\]
Note that $\alpha_{\mathrm{mah}}(s_i)\to 1$ indicates that the structural profile of $s_i$ is statistically typical of some existing area of $\mathcal{L}$, whereas $\alpha_{\mathrm{mah}}(s_i) = \alpha_{\min}$ indicates that $s_i$ is a structural outlier relative to every cluster in $\mathcal{L}$.     
\end{definition}

\subsection{An Analogue of Perplexity for Quantifying Proof Non-triviality}

We now combine these two notions of recoverability to measure whether a machine-generated proof step is truly surprising~[10], or a reformulation of known results.

\begin{definition}[Combined Recoverability]\label{combinedr}
The combined recoverability of a step, denoted as $\alpha(s_i,B)$, is defined as the minimum of its ablation and Mahalanobis recoverability scores.
\[
\alpha(s_i,B) = \min\big(\alpha_{\mathrm{abl}}(s_i,B),\ \alpha_{\mathrm{mah}}(s_i)\big)
\]  
The step's surprisal, $\mathrm{Surp}(s_i,B)$, is then calculated as the negative logarithm of this combined score, bounded by a maximum value $\kappa_{\max} = -\log\alpha_{\min}$.
\[
\mathrm{Surp}(s_i,B) = -\log\alpha(s_i,B)\ \in[0,\kappa_{\max}].
\]
\end{definition}

Taking the minimum over the two independent estimators means a step is only scored as \emph{recoverable} when both the ablation and the distributional (Mahalanobis) signals agree; either estimate alone can cancel recoverability. Having quantified the surprisal of individual steps, we now aggregate these values across the entire proof sequence to establish our formal analogue of perplexity.

\begin{definition}[Proof Perplexity]\label{proofp}
Proof perplexity, denoted as $\mathrm{PP}_w(\pi,B)$, aggregates the surprisal across the entire proof $\pi$. It is computed as the exponential of the weighted average of the individual step surprisals, where each step's weight $w_i=w(s_i)$ corresponds to its elaborated kernel-term size. More precisely,
\[
\mathrm{PP}_w(\pi,B) = \exp\left(\frac{\sum_i w_i\mathrm{Surp}(s_i,B)}{\sum_i w_i}\right).
\]
\end{definition}

The proof perplexity gives us an aggregated score for the entire proof. However, to directly answer our initial question of how much of the machine-generated proof is a novel discovery, we need to measure the exact proportion of the proof that is highly surprising. For this, we set a recoverability threshold $\theta$ and define the following. 

\begin{definition}[$\theta$-invention fraction]\label{ifraction}
For a given recoverability threshold $\theta \in (\alpha_{\min}, 1)$, the $\theta$-invention fraction, $\beta(\theta,B)$, measures the proportion of the proof's total weight that consists of highly surprising steps. It is calculated by summing the weights of all steps that have a combined recoverability score less than or equal to $\theta$, divided by the total weight of all steps in the proof. More precisely,
\[
\beta(\theta,B) = \frac{\sum_{i:\alpha(s_i,B)\leq\theta} w_i}{\sum_i w_i}. 
\]
\end{definition}

The following result shows that if a proof has a high overall perplexity, it guarantees a minimum lower bound on the fraction of the proof that is highly surprising, answering our initial question.

\begin{theorem}\label{thm2.1}
Let $\kappa(\theta):=-\log\theta$. Whenever $\log\mathrm{PP}_w(\pi,B) > \kappa(\theta)$,
\[
\beta(\theta,B)\geq 
\frac{\log \mathrm{PP}_w(\pi,B) - \kappa(\theta)}{\kappa_{\max}-\kappa(\theta)}. 
\] 
\end{theorem}

The proof of the theorem is given in Appendix~\ref{theoremproof}. 

\subsection{The Experimental Result}\label{section3}
We set $B = 50$, $|\mathrm{Cand}_B|=500$, giving $\alpha_{\min}\approx 0.0020$, $\kappa_{\max}\approx 6.2166$, and $\theta = 0.2$. The topological order $\pi$ contains all $N = 2399$ kernel declarations in the dependency closure of the zeta23 proof\footnote{Mathlib4 git repo: \url{https://github.com/leanprover-community/mathlib4}}. Since $\pi$ orders the full closure, every local dependency tree $\mathcal{T}(s_i\mid P_i)$ has depth one and $w(s_i)=|\mathrm{Cite}(s_i)|$ is the elaborated kernel-term size (see Section~\ref{a3}). The Mahalanobis calibration uses a multivariate $t$-null with per-cluster degrees of freedom $\hat{\nu}_c$ and intersection-union aggregation ($\max_c p_c$, Definition~\ref{mahr}). Scoring against the current \texttt{MATHLIB} git repository yields the results given in Table~\ref{onlytable}.

\begin{table}[ht]\label{onlytable}
\caption{Aggregate scores for the zeta23 proof under $B = 50$, $|\mathrm{Cand}_B| = 500$, $\theta = 0.2$.}
\centering
\begin{tabular}{ll}
\toprule
Metric & Value \\
\midrule
Proof perplexity $\mathrm{PP}_w(\pi, B)$ & 388.23 (cap: 501) \\
$\theta$-invention fraction $\beta(\theta, B)$ & 0.9578 (bound: 0.9446, holds) \\
Steps with $\alpha \leq \theta$ / total & 2032\,/\,2399 \\
\bottomrule
\end{tabular}
\end{table}

Of the 2399 declarations, 1992 sit at the floor $\alpha_{\min}$. Only 367 declarations come back with $\alpha>0.2$. Table~\ref{tab21} isolates five non-floor steps that illustrate the benchmark's discriminative behaviour. 

\begin{table}[ht]\label{tab21}
\caption{Representative high-surprisal declarations from the zeta23 kernel closure. Combined recoverability $\alpha \leq 0.2$ is the highly-surprising cut. The $\dagger$~marks a classical result formalised for the first time within the zeta23 project.}
\centering
\begin{tabular}{lrrlllr}
\toprule
Proof step & $w$ & $w_B$ & $\alpha_{\mathrm{abl}}$ & $\alpha_{\mathrm{mah}}$ & $\alpha$ & Surp \\
\midrule
hasDerivAt\_Zeta0Integral & 64\,703\,611 & 15 & 1.0000 & 0.0023 & 0.0023 & 6.0782 \\
Zeta23.EF.pole\_term & 886 & 24 & 0.9729 & 0.0032 & 0.0032 & 5.7498 \\
Zeta23.Cheb.sum\_vonMangoldt\_sq\_le & 38 & 5 & 0.8684 & 0.0115 & 0.0115 & 4.4630 \\
RHLinalg.rank\_trace\_ineq & 150 & 13 & 0.9133 & 0.0480 & 0.0480 & 3.0362 \\
RHLinalg.vonNeumann\_trace\_ineq$^{\dagger}$ & 85 & 85 & 0.0020 & 0.4457 & 0.0020 & 6.2166 \\
\bottomrule
\end{tabular}
\end{table}

The first four rows share a common signature. Ablation search partially or fully recovers the kernel term ($\alpha_{\mathrm{abl}}$ is high), yet the Mahalanobis test places the step far from every cluster in $\mathcal{L}$ ($\alpha_{\mathrm{mah}}$ is small). Combined recoverability is driven by the structural outlier signal, and $\alpha$ remains well below $\theta = 0.2$. This is precisely the pattern the combined minimum rule was designed to detect, computationally accessible material whose mathematical argument is nonetheless structurally alien to the existing library.

RHLinalg.rank\_trace\_ineq is the cleanest example. The prover rebuilds 13 of 150 internal units ($\alpha_{\mathrm{abl}} = 0.9133$), confirming that the step's ingredients are individually present in \texttt{MATHLIB}. Yet the Mahalanobis test returns $\alpha_{\mathrm{mah}} = 0.0480$. The way those ingredients are combined, a rank-trace inequality that bounds the number of large eigenvalues of a Hermitian perturbation has no statistical precedent in any cluster of $\mathcal{L}$. Anthropic's own technical report~[4] singles out this inequality as the novel analytic input of the two-thirds argument, and the benchmark independently confirms that assessment. 


RHLinalg.vonNeumann\_trace\_ineq illustrates a limitation inherent in any \texttt{MATHLIB} relative benchmark. Von Neumann's trace inequality is a classical result (1937), but it was not formalised in \texttt{MATHLIB} before the zeta23 project. Since it is absent from $\mathcal{L}$, ablation search cannot find it ($w(s_i, B) = w(s_i)$, $\alpha_{\mathrm{abl}} = \alpha_{\min}$), and the step registers at the surprisal ceiling. Note, however, that $\alpha_{\mathrm{mah}} = 0.4457$, the Mahalanobis test correctly places it within the statistical neighbourhood of existing linear-algebra declarations, confirming that its structural profile is not novel even though its formalisation is. The combined minimum still pins it at the floor because ablation is binding. This is a deliberate conservatism: $\beta(\theta, B)$ should be read as a certificate of non-recoverability from \texttt{MATHLIB}, not as a count of new mathematical ideas. Results like this inflate the invention fraction, and we expect the effect to diminish as \texttt{MATHLIB} grows.

\newpage
\section*{References}

{
\small
[1] Somani, N.\ (2026) PriorProof: A Point-in-Time Measure of Technique Novelty for Formal Proofs, {\it Pre-print: arXiv:2607.16997v1}, pp.\ 1--17.

[2] Constantine, G., Pan, L., \& Ganesh, V. (2026). Market Incentivization for Theorem Proving with LLM Agents. {\it Proceedings of the 43rd International Conference on Machine Learning (ICML 2026)}, PMLR 306. 

[3] Wang, Z., Dai, Q., Wu, Y., \& Wen, Z. (2026). FaithSieve: Fine-Grained Evaluation of Math Proofs with Faithful Formal Evidence. {\it Pre-print: arXiv:2608.26310}, pp. \ 1 -- 32.

[4] Alp\"{o}ge, L., \& Furman, R. (2026). More than two thirds of the zeros of the Riemann zeta function are simple and on the critical line. {\it Pre-print: arXiv:2608.13637v1}, pp. \ 1 -- 17.

[5] Miaschi, A.; Brunato, D., Dell'Orletta, F. \& Venturi, G. (2021). ``What Makes My Model Perplexed? A Linguistic Investigation on Neural Language Models Perplexity". {\it Proceedings of Deep Learning Inside Out (DeeLIO): The 2nd Workshop on Knowledge Extraction and Integration for Deep Learning Architectures.} pp. 40 -- 47.

[6] Ledoit, O., \& Wolf, M. (2004). A well-conditioned estimator for large-dimensional covariance matrices. {\it Journal of Multivariate Analysis}, {88(2)}, pp. \ 365 -- 411.

[7] Bonferroni, C. (1936). Teoria statistica delle classi e calcolo delle probabilit\'{a}. {\it Pubblicazioni del R. Istituto Superiore di Scienze Economiche e Commerciali di Firenze}, 8, pp. \ 3 -- 62.

[8] Hubert, T. et al (2026). Olympiad-level formal mathematical reasoning with reinforcement learning. {\it Nature} 651, pp. \ 607 -- 613.

[9] Achim, T., Best, A., Bietti, A., et al. (2025). Aristotle: IMO-level Automated Theorem Proving. {\it Pre-print: arXiv:2510.01346}, pp \ 1 -- 23.

[10] Futrell, R., Gibson, E., \& Levy, R. P. (2020). Lossy-Context Surprisal: An Information-Theoretic Model of Memory Effects in Sentence Processing. {\it Cognitive Science}, {44}, pp. \ 1 -- 54.

[11] Li, X., Peng, N., Severini, S., \& Shafto, P. (2026). The network structure of Mathlib. {\it Pre-print: arXiv:2604.24797}, pp. \ 1 -- 66.

[12] Mahalanobis, P. C. (1936). On the generalized distance in statistics. {\it Proceedings of the National Institute of Sciences of India}, (2), pp. \ 49 -- 55.

[13] Otter, R. (1949). The multiplicative process. {\it Ann. Math. Statist.,} 20(2), pp. \ 206 -- 224.

[14] Zipf, G. K. (1949). Human Behavior and the Principle of Least Effort. {\it Addison-Wesley.}

[15] Anthropic. (2026). formal-math: Machine-checked Lean 4 formalizations published by Anthropic. GitHub. {\url{https://github.com/anthropics/formal-math}}

[16] Kotz, S., \& Nadarajah, S. (2004). Multivariate $t$ Distributions and Their Applications. Cambridge University Press.
}


\appendix

\section{Further Details on the Methodology}

\subsection{Note on the Related Work} One of the closest related approaches was developed by N. Somani, known as Prior Proof~[1]. This approach measures proof route novelty relative to the historical state of the \texttt{MATHLIB} library. For a theorem introduced at time $t$, it constructs the prior from the \texttt{MATHLIB} snapshot preceding the start of the quarter containing $t$, and scores the proof according to how surprising its dependency families are under that historical prior. Thus, the author studies the novelty of a proof step explicitly as a time-relative object. A proof route can be novel when introduced but become less novel as the same dependencies become more prevalent in later proofs.

While Prior Proof captures how the novelty of a specific proof route changes over time, its reliance on historical point-in-time snapshots limits its ability to measure the structural uniqueness of a proof step or how easily it can be rediscovered by advanced automated search algorithms. Our approach differs by moving away from temporal baselines and instead evaluating proof steps through a probabilistic and reconstruction perspective. We only mention the above result in detail as it is closely related to our work but differs in the foundational methodology.

\subsection{Proof Tree Structure in the Mathlib Library}

\subsubsection{Lean Kernel Declarations and Dependencies}

In this section, we make precise what a Lean declaration is and what it means for one declaration to depend on another, since these notions underlie both $\mathcal{L}$ and the topological order $\pi$ introduced in Section~\ref{1.1}.

Let $\Gamma$ denote the ambient Lean environment relevant to our setting, i.e. $\Gamma = \mathcal{L}\cup \{T,s_1, \ldots, s_n\}$: the fixed \texttt{MATHLIB} library together with the target theorem $T$ and its intermediate proof steps. We take $\Gamma$ to be finite in  \texttt{MATHLIB}, together with any single theorem's addition to it, is a concrete, finite collection of declarations at any point in time. A standing hypothesis used below to guarantee termination of the recursive constructions in Definitions~\ref{localdtree} and~\ref{searchtreeformal}.

\begin{definition}[Lean Declaration]
A Lean declaration $d \in \Gamma$ is a triple $d = (n_d,\ \tau_d,\ e_d)$,  where $n_d$ is the unique name of $d$ in $\Gamma$, $\tau_d$ is its type for a lemma, proposition, or theorem, $\tau_d$ is the $\mathrm{Prop}$-valued statement being asserted, and $e_d$ is the fully elaborated term inhabiting $\tau_d$, i.e. the proof term itself. The Lean kernel accepts $d$ into $\Gamma$ only after independently type-checking that $e_d:\tau_d$.
\end{definition}

Kernel checking of $e_d : \tau_d$ proceeds by elaborating every constant symbol that occurs in $e_d$ or $\tau_d$. Each such constant must itself already be a checked declaration in $\Gamma$; the kernel does not re-derive it, it looks it up. This lookup relation is exactly what we mean by dependency.

\begin{definition}[Direct Dependency]\label{directdep}
For $d, d' \in \Gamma$, we say $d'$ is a direct dependency of $d$, written $d' \prec d$, if $n_{d'}$ occurs as a constant symbol in $e_d$ or $\tau_d$. Equivalently, the kernel must have already admitted $d'$ into $\Gamma$ before $d$ can be type-checked.
\end{definition}

We write $\prec^{}$ for the irreflexive transitive closure of $\prec$: $d' \prec^{} d$ holds iff there is a non-empty chain $d' \prec d^{(1)} \prec \cdots \prec d^{(k)} \prec d$ for some $k \ge 0$, with $k=0$ meaning $d' \prec d$ directly. In particular, $d \prec^{*} d$ would require a $\prec$-cycle through $d$.

\begin{lemma}[$(\Gamma, \prec)$ is a DAG]\label{dag}
The relation $\prec$ admits no cycles: $(\Gamma, \prec)$ is a directed acyclic graph (DAG)\footnote{For a proof of this lemma, see appnexdix A of~[11]}.    
\end{lemma}

Lean's elaborator processes declarations as an ordered sequence and forbids forward references. A declaration may only cite names already present in $\Gamma$ at the point it is elaborated recursive declarations are handled separately via explicit well-founded recursion principles, which are themselves ordinary prior declarations. Hence there is no $d$ with $d \prec^{*} d$.

\begin{definition}[Dependency Closure]\label{depclose}
For $d \in \Gamma$, the dependency closure of $d$ is
\[
\mathrm{Dep}(d) := \{, d' \in \Gamma : d' \prec^{*} d ,\},
\]
the set of all declarations the kernel must have checked, directly or transitively, in order to accept $d$. 
\end{definition}

By Lemma~\ref{dag}, $d \notin \mathrm{Dep}(d)$. This is precisely the closure referenced in Section~\ref{1.1}. For the target theorem $T$, $\mathrm{Dep}(T)$ is ``the set of Lean declarations on which the proof steps of $T$, as checked by the Lean kernel, actually depend.'' Note that $\mathrm{Dep}(d)$ is unfolded in full, all the way down to $\Gamma$'s foundational declarations; Below we introduce a relative version of this closure, truncated at $L \cup \{s_1,\ldots,s_{i-1}\}$, which is what $w(s_i)$ is actually built from.

{\it Recovering the Topological order $\pi$.} Because $(\Gamma, \prec)$ is acyclic, any linear extension of $\prec$ restricted to $\{s_1, \ldots, s_n\}$ is well-defined. The order $\pi = (s_1, \ldots, s_n)$ of Section~\ref{1.1}, satisfying $\mathrm{Dep}(s_i) \cap \{s_1,\ldots,s_n\} \subseteq \{s_1,\ldots,s_{i-1}\}$, is exactly such a linear extension; it is guaranteed to exist by Lemma~\ref{dag} and is generally not unique.

\subsection{Lean Kernel Declarations and the Internal Complexity Tree}\label{a3}

We now unfold the DAG $(\Gamma, \prec)$ into a tree in order to give $w(s_i)$ and $w(s_i, B)$ a precise combinatorial meaning.

\begin{definition}[Local Dependency Tree]\label{localdtree}
Fix a step index $i$ and let $P_i := \mathcal{L} \cup \{s_1, \ldots, s_{i-1}\}$ denote the declarations already available prior to $s_i$ in the topological order $\pi$ (so that $P_1 := \mathcal{L}$).
\end{definition}

For $d \in \Gamma$, let $\mathrm{Cite}(d) = (c_1, c_2, \ldots, c_m)$ be the sequence of constant occurrences in $e_d$ and $\tau_d$, read in a fixed left-to-right pre-order traversal of the elaborated term's syntax, where each $c_j \in \Gamma$ is the declaration referenced at that occurrence. This is deliberately different from the direct-dependency relation $\prec$ of Definition~\ref{directdep}, which only records whether $d' \prec d$ holds and so collapses repeated citations into a single fact. Observe that $\mathrm{Cite}(d)$ is a sequence, not a set. If the same declaration is invoked at several distinct positions in $e_d$, it appears several times in $\mathrm{Cite}(d)$, once per occurrence. The two notions are related by
\[
d' \prec d \iff d' \mathrm{ occurs\; at\; least\; once\; in }\;\mathrm{Cite}(d),
\]
i.e. $\prec$ is the multiplicity-forgetting shadow of $\mathrm{Cite}$; it remains the right notion for Lemma~\ref{dag} and $\mathrm{Dep}(d)$ (see Definition~\ref{depclose}), where only reachability, not cost, matters.

Define the local dependency tree of $d$ relative to $P_i$, written $\mathcal{T}(d \mid P_i)$, recursively:
\begin{itemize}
\item the root is labeled $d$;
\item if $d \in P_i$, the tree is just this root, with no branches, a declaration already present in the library is treated as atomic and is not expanded further;
\item otherwise, writing $\mathrm{Cite}(d) = (c_1, \ldots, c_m)$, the root's branches are the trees $\mathcal{T}(c_1 \mid P_i), \ldots, \mathcal{T}(c_m \mid P_i)$, one branch per occurrence, in order. If some declaration $c$ is cited $k$ times in $e_d$, it contributes $k$ separate (isomorphic) branch subtrees, not one.
\end{itemize}
Since each $e_d$ is a finite syntactic term (so each $\mathrm{Cite}(d)$ is a finite sequence) and $(\Gamma, \prec)$ is acyclic, therefore, no chain of non-$P_i$ expansions can recur, this recursion terminates and $\mathcal{T}(d \mid P_i)$ is a finite tree.

A declaration $c \notin P_i$ can label more than one node of $\mathcal{T}(s_i \mid P_i)$ for two different reasons, and both are counted: 
(i) Occurrence multiplicity $c$ is cited more than once within a single ancestor's term, contributing repeated siblings by construction above; and 
(ii) sharing multiplicity $c$ is cited by two or more distinct ancestors along different branches of the tree (the usual DAG-to-tree unfolding phenomenon). Both represent genuine, separately incurred proof effort citing a fact twice costs two applications, whether those two applications sit under the same parent step or under two different ones, so neither is collapsed.

Unlike $\mathrm{Dep}(d)$, which unfolds all the way down to $\Gamma$'s foundational declarations regardless of what is already established, $\mathcal{T}(d \mid P_i)$ stops as soon as it reaches material already available in $P_i$. This is not a simplification for convenience, it is necessary. Every proof step $s_i$ necessarily cites a large amount of pre-existing library material (such that arithmetic identities, order lemmas, etc.); fully unfolding that material, as the naive closure $\mathrm{Dep}(s_i)$ would, swamps $w(s_i)$ with complexity that is common to essentially every theorem in \texttt{MATHLIB} rather than specific to $s_i$. Truncating at $P_i$ isolates exactly the declarations genuinely internal and new to $s_i$'s own argument matching our description of $w(s_i)$ as ``the internal declarations of the proof step $s_i$".

Formally, we write the internal complexity of a proof step $s_i$ is the number of non-root nodes of its local dependency tree:
\[
w(s_i) := \big|\mathrm{Nodes}\big(\mathcal{T}(s_i \mid P_i)\big)\big| - 1.
\]

By Definition~\ref{localdtree}, this counts, with multiplicity, once per syntactic occurrence in $e_{s_i}$, and once per distinct sharing branch, every citation reachable from the original, already-formalised proof term $e_{s_i}$. A citation of a declaration already in $P_i$ contributes exactly one unit each time it is invoked (it is not itself expanded), while a citation of a genuinely new auxiliary declaration introduced specifically within $s_i$'s proof is expanded recursively and contributes its own internal structure, again once per occurrence. In particular, $w(s_i)$ is bounded by the true size of $e_{s_i}$'s novel content, not by the size of \texttt{MATHLIB} beneath it.

When $\pi$ orders the entire dependency closure of $T$, as in our experiment, where $N = 2399$, every declaration that $s_i$ cites already lies in $P_i$ by topological ordering. Consequently, every child of the root in $\mathcal{T}(s_i \mid P_i)$ is a $P_i$-leaf, the recursive expansion terminates after a single level and no child of the root requires further expansion. Additionally, $\mathcal{T}(s_i \mid P_i)$ is a depth-one fan-out. In this regime $w(s_i) = |\mathrm{Cite}(s_i)|$ is the number of constant occurrences in the elaborated term $e_{s_i}$, which is precisely the elaborated kernel-term size referenced in Definition~\ref{proofp}. The recursive generality of Definition~\ref{localdtree} remains available for experimental designs in which $\pi$ contains only a proper subset of the closure (for instance the 33 top-level named lemmas), so that intermediate declarations are expanded as internal tree nodes; in the present experiment, however, $w(s_i)$ should be read as a flat citation count rather than a recursive tree size.

The search tree $\mathcal{T}(s_i)$ of Definition~\ref{searchb} is a different object from $\mathcal{T}(s_i \mid P_i)$: one is the ground-truth structure of an already-existing proof, the other is a stochastic exploration by an automated prover attempting to rediscover a proof of $s_i$ from scratch. But the two must be and, by construction, are counted in the same units. Both charge one unit per action, occurrence or application, never deduplicating repeated uses of the same fact. We make this precise below.

\begin{definition}[Search Tree Structure]\label{searchtreeformal}
The search tree $T\mathcal{T}s_i)$ is a rooted tree whose nodes are proof states and whose edges are tactic applications. A proof state $\sigma = (g_1, \ldots, g_m)$ is an ordered list of currently open goals. The root is the singleton state $(s_i)$. Each edge leaving a state $\sigma$ corresponds to selecting one open goal $g_j \in \sigma$ and applying a single tactic or lemma drawn from $P_i$, filtered in practice to the curated candidate pool $\mathrm{Cand}_B \cup \{s_1, \ldots, s_{i-1}\}$ for tractability. If the applied lemma has $k$ premises, $g_j$ is replaced by $k$ new subgoals in the successor state; if $k = 0$, $g_j$ is simply removed. Branching at a node represents alternative tactic choices for the same proof state (OR-branching); the generation of multiple subgoals from a single application is encoded within the successor state, not as separate branches. Crucially, applying a declaration from $P_i$ costs exactly one edge regardless of that declaration's own internal proof size. The prover never re-derives it, mirroring the leaf convention of the local dependency tree. A root-to-node path $p$ has length $|p|$ equal to its number of edges, i.e. the total number of tactic applications across all goals, and the search is truncated at depth $B$. Any branch with $|p|>B$ is discarded, so every explored path satisfies $|p|\leq B$. A path $p$ is accepted if its terminal state is the empty list $()$, meaning every goal has been closed, yielding a kernel checkable proof of $s_i$ from $P_i$ alone. 
\end{definition}

{\it Recovered Complexity $w(s_i, B)$.} If $\mathcal{T}(s_i)$ contains at least one accepted path within depth $B$, let $p^{\ast}$ be the accepted path found by the prover, and set $w(s_i, B) := |p^*|$, the total number of tactic applications in this newly recovered proof. If no accepted path exists within the budget then the search tree is exhausted to depth $B$ without closing the goal. In such cases, the prover has failed to rediscover $s_i$ under budget $B$, and we set, as in Section~\ref{1.1},
$$w(s_i, B) := w(s_i).$$

Setting $w(s_i,B) := w(s_i)$ on failure is not an approximation but a deliberate convention: it forces the ratio $w(s_i,B)/w(s_i) = 1$ in Definition~\ref{ablationr}, so $\alpha_{\mathrm{abl}}(s_i,B)$ collapses to its floor $\alpha_{\min}$. This is the intended behaviour because a step the automated prover cannot rediscover at all within budget $B$ should register as maximally non-recoverable, on the same footing as a step whose relative recovered size is already at the floor.

{\it Relating the Two Tree Structures.} $\mathcal{T}(s_i \mid P_i)$ and $\mathcal{T}(s_i)$ coincide only in the degenerate case where the automated prover's accepted path $p^{\ast}$ reconstructs a proof isomorphic to $s_i$'s original one. In general $\mathcal{T}(s_i \mid P_i)$ is fixed and it reflects what \texttt{MATHLIB} library already contains, while $\mathcal{T}(s_i)$ is stochastic and budget-bounded. What makes them comparable is that both are relativised to the same boundary $P_i$ and both charge one unit per proof action. Observe that $w(s_i)$ counts one non-root node per citation occurrence in the dependency tree, while $|p^*|$ counts one edge per tactic application in the proof-state tree, so $w(s_i)$ and $w(s_i, B)$ are measured in comparable units.

\subsection{Reflection on the Mahalanobis Recoverability}

In this section, we first write the motivation for using Mahalanobis distance~[12], and then we provide further details on the methodology.

\subsubsection{Motivation for using Mahalanobis Distance} \label{appendixformah}
For a scalar feature $X$ with mean $\mu$ and standard deviation $\sigma$, recall the standardised third moment, or skewness, $\gamma_1 := \mathbb{E}[(X-\mu)^3]/\sigma^3$. Symmetric distributions, and in particular the Gaussian, have $\gamma_1 = 0$.

Within a fixed mathematical topic, however, we should not expect proof-step complexity to be distributed symmetrically around its mean. Most declarations are routine: they are short, standard applications of a small number of familiar techniques. A smaller number are more exceptional, combining several non-trivial prior results into something structurally deeper. The resulting population is therefore right-skewed (i.e., $\gamma_1>0$), with a dense cluster of low-complexity, routine declarations and a long right tail of unusually deep ones. In other words, we should expect something quite different from a bell curve centred symmetrically around some ``typical'' level of complexity.

This asymmetry is exactly why raw Euclidean distance $\|\phi(s_i)-\mu_c\|$ is a poor measure of typicality. Euclidean distance treats a deviation of magnitude $\Delta$ above the mean as equally unusual as a deviation of the same magnitude below it. But in a skewed population, that symmetry does not hold. A step $\Delta$ above $\mu_c$ may be completely unremarkable, lying well within the population's long right tail, while the same step $\Delta$ below $\mu_c$, where there is much less room for variation because complexity is bounded below near zero, may already be genuinely rare. A distance measure that ignores the shape of the population it is applied to can therefore misjudge both the direction and the magnitude of what should actually count as surprising.

Mahalanobis distance replaces this raw deviation with a whitened one,
\[
D_M(s_i \mid c) = \sqrt{(\phi(s_i)-\mu_c)^{\top}\Sigma_c^{-1}(\phi(s_i)-\mu_c)},
\]
as given in Section~\ref{2.1.1}, where $\Sigma_c$ is fit directly from $\{\phi(d) : d \in L,, d \in c\}$ rather than assumed in advance. Because $\Sigma_c$ is estimated from the actual population, it reflects the empirical spread of cluster $c$ along each feature axis, including the elongated spread produced by a heavy right tail. Concretely, if the complexity axis has a particularly heavy tail within cluster $c$, the corresponding entry of $\Sigma_c$ will also be large, so $\Sigma_c^{-1}$ down weights deviations along that axis. This is the correction Euclidean distance is missing, a step that falls within the long but ordinary tail of a heavy tailed cluster is not penalised as though it were just as unusual as an equal magnitude deviation along an axis where the population is tightly concentrated.

We use Ledoit-Wolf shrinkage because this mathematically stabilises the covariance matrix $\Sigma_c$ to ensure it is invertible (which is strictly required to calculate $\Sigma_c^{-1}$ in the Mahalanobis distance). In high dimensional data like proof embeddings, standard covariance matrices often break or become unstable. Ledoit-Wolf shrinkage fixes this by blending the empirical data with a structured baseline.


Under a multivariate $t$ with $\nu_c$ degrees of freedom~[\S2, 16],
the quadratic form built from the true scale matrix follows
$q\cdot F_{q,\nu_c}$, whose tail decays polynomially as
$P(F_{q,\,\nu_c} > t) \sim t^{-\nu_c/2}$, matching the structural
prediction. Because the Ledoit-Wolf estimator~[6] returns
the covariance matrix rather than the scale matrix, and the two differ by the factor $\nu_c / (\nu_c - 2)$ for $\nu_c > 2$, the Mahalanobis distance computed from the covariance estimate must be rescaled accordingly
\[
\frac{\nu_c}{q(\nu_c - 2)}\,D_M(s_i\mid c)^2 \;\sim\; F_{q,\,\nu_c}\,,
\]
giving the tail probability
\[
p_c(s_i) = 1 - F_{q,\,\nu_c}\!\left(\frac{\nu_c}{q(\nu_c - 2)}
\,D_M(s_i\mid c)^2\right)
\]
as stated earlier. 

The parameter $\nu_c$ is estimated per cluster from the empirical excess kurtosis of within cluster $D_M(d \mid c)^2$ values, by inverting the theoretical kurtosis of the $F_{q,\nu}$ distribution.

We aggregate across clusters using an intersection–union test (IUT). The relevant question is whether $s_i$ is typical of no cluster, i.e., whether it is a structural outlier relative to every topic area in $\mathcal{L}$. The IUT rejects typicality when $\max_c p_c(s_i) < \delta$. Under the global null $H_0:\exists\, c_0$ such that $s_i$ is a typical draw from cluster $c_0$, the inequality $\max_c p_c \geq p_{c_0}$ gives $P_{H_0}(\max_c p_c<\delta) \leq P_{H_{0,c_0}}(p_{c_0} <\delta) \leq \delta$, since $p_{c_0}$ is stochastically at least $\texttt{Uniform}[0,1]$ under the true null $H_{0,c_0}$. This controls the type I error at level $\delta$ without multiplicity corrections such as Bonferroni~[7].



\subsubsection{Limitations}

Note that $\Sigma_c$ is a second-moment object. It captures the scale of the tail in each direction, but not its asymmetry. The multivariate~$t$ null used to calibrate $D_M(s_i\mid c)^2$ via $F_{q,\,\nu_c}$ still assumes an elliptically symmetric reference population. Thus, if a cluster is genuinely skewed, whitening does not make it exactly $t$-distributed, and the resulting tail probability $p_c(s_i)$ should be understood as an approximation rather than an exact probability. We treat this as a deliberate and acknowledged simplification, rather than an oversight. It is precisely the kind of higher-order residual that Ledoit-Wolf shrinkage and the intersection-union aggregation are not intended to remove.

\subsubsection{Mahalanobis Distance and Lean kernel Structure}\label{a43}

The skewness described above is not simply an incidental empirical feature of the mathematical metadata in \texttt{MATHLIB}. It has a structural origin in the dependency and citation patterns of the Lean kernel.

In our experiment, where $\pi$ orders the full dependency closure, every local dependency tree $\mathcal{T}(s_i \mid P_i)$ has depth one, so $w(s_i)=|\mathrm{Cite}(s_i)|$ is the constant occurrence count of the elaborated term. Even in this flat regime the distribution of $w(s_i)$ within a cluster is right skewed for a concrete reason. Elaborated kernel terms vary enormously in size, from single line definitions to large analytic estimates, and the population within any topic area is dominated by a bulk of short routine declarations with a long right tail of substantially larger ones\footnote{When $\pi$ contains a proper subset of the closure and Definition~\ref{localdtree} recurses beyond one level, classical branching process theory~[13] predicts that multiplicative branching further amplifies this right skewness; in the full-closure regime of our experiment, however, all trees have depth one and this mechanism is inert.}. In the zeta23 closure $w(s_i)$ ranges from 2 to 68,703.

The premise usage side of the feature vector reinforces this asymmetry. Dependency and citation structures in formal libraries exhibit a Zipf-like [14] usage pattern. A small number of foundational declarations (arithmetic identities, order lemmas, basic topology facts) are cited by a large fraction of the library, while many specialised declarations are cited only rarely. The resulting premise usage histograms $\phi^*_{\mathrm{usage}}(s_i)$ are therefore sparse and heavy tailed across the library as a whole, and right skewed within individual topic clusters.

Thus Mahalanobis whitening is not being applied to an arbitrary numerical summary. It is being applied to quantities whose skewness is a direct consequence of the heterogeneous elaborated term sizes and the Zipf-like citation structure of the Lean kernel.

{\it Structural Feature Vector of the Proof Steps.} Section~\ref{2.1.1} allows $\phi(s_i)$ to be constructed a specific way, and this construction comes directly by using the machinery of the Lean's kernel structure~[11].

Fix an enumeration of $P_i$ or, for tractability, of $\mathrm{Cand}_B$. Define $\phi*{\mathrm{usage}}(s_i) \in \mathbb{N}^{|P_i|}$ so that the coordinate corresponding to $d$ counts how many times $d$ appears as a $P_i$-leaf in $\mathcal{T}(s_i \mid P_i)$. In other words, this is a leaf-multiplicity histogram that can be read directly from the occurrence sequence $\mathrm{Cite}$ as it is expanded through the tree.

A coarser, unweighted representation simply drops these counts. Specifically, $(\phi_{\mathrm{dep}}(s_i))_d = {\bf 1}[d \in \mathrm{Dep}(s_i) \cap P_i]$. This vector is obtained from the full dependency closure $\mathrm{Dep}(d)$, restricted to $P_i$.

In the full closure regime $w(s_i) = |\mathrm{Cite}(s_i)|$ is a scalar, so the natural complexity embedding is $w(s_i)$ itself the elaborated kernel-term size. Concatenating it with $\phi^*_{\mathrm{usage}}(s_i)$ or $\phi_{\mathrm{dep}}(s_i)$ gives a combined vector that captures both the scale of the proof term and its premise usage profile. In experimental designs where $\pi$ contains a proper subset of the closure and $\mathcal{T}(s_i \mid P_i)$ has depth greater than one, $w(s_i)$ can be extended to a small vector of tree-shape statistics (such as depth, maximum branching factor, and the fraction of non-$P_i$ internal nodes), giving a richer embedding of the proof step's structural complexity.


\subsubsection{Motivation behind Partitioning Mathlib into Clusters}

Mahalanobis distance is meaningful only relative to a well-defined reference population $(\mu_c, \Sigma_c)$. Fitting a single global $(\mu, \Sigma)$ across all of $\mathcal{L}$ would combine declarations from structurally different areas, such as elementary arithmetic and deep measure theory, into a single population. This is problematic because the skewness described above is itself a local phenomenon. Each topic area has its own characteristic level of \emph{routine} complexity, as well as its own local heavy tail of unusually deep results.

A global fit would therefore conflate two different kinds of variation. Genuine differences between topics, where results in this area are inherently deeper, and novelty within a topic, where this result is unusually deep even by the standards of its own area. The pooled $\Sigma$ would be inflated by this between cluster heterogeneity, making $D_M$ less sensitive precisely where sensitivity matters most. When a step is a genuine outlier relative to its own mathematical neighbourhood, even if it appears unremarkable when compared with declarations from an unrelated area that is naturally more complex.

This gives a concrete interpretation to the \emph{statistical neighbourhood}. A neighbourhood is meaningful only when it is defined relative to a reasonably homogeneous reference population, which clustering by module or topic provides.

\subsection{Additional Details on the Notion of Proof Perplexity}

\subsubsection{Motivation for Combining the Recoverability Scores
}

Write $\mathrm{Surp}*{\mathrm{abl}}(s_i,B) := -\log \alpha*{\mathrm{abl}}(s_i,B)$ and $\mathrm{Surp}*{\mathrm{mah}}(s_i) := -\log \alpha*{\mathrm{mah}}(s_i)$ for the surprisal each estimator would assign on its own. Since $-\log(\cdot)$ is strictly decreasing, taking the minimum of two recoverability scores is exactly equivalent to taking the maximum of the two surprisals
\[
\mathrm{Surp}(s_i,B) = -\log\big(\min(\alpha_{\mathrm{abl}}(s_i,B), \alpha_{\mathrm{mah}}(s_i))\big) = \max\big(\mathrm{Surp}_{\mathrm{abl}}(s_i,B),\ \mathrm{Surp}_{\mathrm{mah}}(s_i)\big).
\]
This identity makes precise what combining the two estimators by a minimum accomplishes. A step is judged \emph{surprising} as soon as either estimator judges it surprising, and judged recoverable only when both agree it is recoverable. Concretely, it rules out one specific failure mode, a step cannot be certified as trivial (low $\mathrm{Surp}$) just because one proxy happened to be uninformative for it, so long as the other proxy correctly flags it as surprising. A step an under tuned automated prover happens to rediscover cheaply (low $w(s_i,B)/w(s_i)$, hence low $\mathrm{Surp}*{\mathrm{abl}}$) but which is nonetheless a genuine structural outlier in $\mathcal{L}$ (high $\mathrm{Surp}*{\mathrm{mah}}$) is still scored as surprising overall, because the maximum in the identity above is driven by the larger of the two.

This combination rule is also exactly what keeps $\mathrm{Surp}(s_i,B)$ within the given range $[0,\kappa_{\max}]$. Both $\alpha_{\mathrm{abl}}(s_i,B) = \max(1-w(s_i,B)/w(s_i),\alpha_{\min})$ and $\alpha_{\mathrm{mah}}(s_i) = \max(\tilde p(s_i),\alpha_{\min})$ are individually constrained to $[\alpha_{\min},1]$ by their own $\max(\cdot,\alpha_{\min})$ constructions. The first one because $w(s_i,B)/w(s_i)\ge 0$ implies $1-w(s_i,B)/w(s_i)\le 1$, and the second one because $\tilde p(s_i)\in[0,1]$ by construction of the maximum $p$-value across clusters. Since $\alpha_{\mathrm{abl}}$ and $\alpha_{\mathrm{mah}}$ share the same floor $\alpha_{\min}$, their minimum is itself guaranteed to lie in $[\alpha_{\min},1]$. Had the two estimators used different floors, the combined score’s range would only be guaranteed down to the smaller of the two, and a single, common $\kappa_{\max} := -\log\alpha_{\min}$ would not cleanly bound $\mathrm{Surp}$ arising from either kind of evidence. Sharing $\alpha_{\min}$ is therefore not just a matter of tidiness; it is what makes one common $\kappa_{\max}$ meaningful for both signals at once, which Theorem~\ref{thm2.1} depends on directly.

\subsubsection{Limitations of Combined Recoverability}

The aforementioned identity makes the combination rule's failure mode precise, by symmetry with what it prevents. Taking a maximum protects against either estimator unilaterally certifying a step as trivial, but provides no corresponding protection against either estimator unilaterally certifying a step as surprising. If ablation search fails for reasons unrelated to genuine mathematical depth such as an under filtered candidate pool, an unlucky search order, or a budget $B$ too small for an otherwise routine step, the resulting inflated $\mathrm{Surp}*{\mathrm{abl}}(s_i,B)$ propagates directly into $\mathrm{Surp}(s_i,B)$, regardless of what $\mathrm{Surp}*{\mathrm{mah}}(s_i)$ says, since the maximum simply takes whichever value is larger. The same holds in reverse for a spuriously high $\mathrm{Surp}_{\mathrm{mah}}$, for example from a small or poorly separated cluster.

This is a genuine asymmetry. A symmetric alternative such as $\alpha(s_i,B) := \sqrt{\alpha_{\mathrm{abl}}(s_i,B),\alpha_{\mathrm{mah}}(s_i)}$, averaging rather than maximising the two surprisals would dampen the effect of a single noisy estimator, but at the cost of losing the ``surprising unless both agree it is not" guarantee that the above identity relies on. A step could then be scored as recoverable simply because one very low surprisal outweighs one moderately high one on average, exactly the outcome the minimum rule was chosen to prevent. Our choice of $\min$ is therefore a deliberate asymmetric trade-off. It favours never dismissing genuine novelty over never overstating it. At the same time, it inherits, rather than removes, whatever miscalibration either individual estimator carries on its own.

\subsubsection{Proof Perplexity from the Surprisal Function}\label{a.5.3}

Taking logarithm in Definition~\ref{proofp},
\[
\log\mathrm{PP}_w(\pi,B) = \frac{\sum_i w_i\,\mathrm{Surp}(s_i,B)}{\sum_i w_i},
\]
so Proof Perplexity is, by construction, the exponential of a $w_i$-weighted average surprisal, exactly the same construction as ordinary perplexity in language modelling, where perplexity is the exponential of the average per-token cross-entropy. Here, a proof step $s_i$ plays the role of a token, while uniform weighting is replaced by $w_i = w(s_i)$, the elaborated kernel-term size, so larger, more substantial steps contribute proportionally more to the aggregate.

This connection is exact at the level of the scale over which $\mathrm{PP}_w$ ranges, not just an evocative label. Recall $\alpha_{\min} = 1/(1+|\mathrm{Cand}_B|)$, it is precisely the probability of selecting, by uniform chance, whichever one of the $|\mathrm{Cand}_B|$ pool candidates together with one null or trivial action actually closes a goal. Consequently, $\kappa_{\max}=-\log\alpha_{\min}=\log(1+|\mathrm{Cand}_B|)$ is exactly the Shannon entropy, in nats, of this same uniform baseline, and $\mathrm{PP}_w(\pi,B) \in \big[\exp(0),\ \exp(\kappa_{\max})\big] = \big[1,\ 1+|\mathrm{Cand}_B|\big]$.

This mirrors the classical reading of perplexity quite closely. A perplexity of $1$ means that every step was recovered with complete confidence ($\mathrm{Surp}=0$ throughout), while a perplexity of $1+|\mathrm{Cand}_B|$ means that, on weighted average, the proof is as surprising as simply guessing uniformly from the candidate pool. In other words, it is the proof-theoretic analogue of a language model whose predictions are no better than picking uniformly from a vocabulary of size $1+|\mathrm{Cand}_B|$.

\subsubsection{The bound of Invention Fraction using Proof Perplexity}

Writing $W := \sum_i w_i$ and treating $w_i/W$ as a probability weight over steps $i=1,\dots,n$, Definition~\ref{ifraction} is exactly the weighted survival function of $\mathrm{Surp}(s_i,B)$ evaluated at $\kappa(\theta)$,
\[
\beta(\theta,B) = \Pr\nolimits_w\big[\alpha(s_i,B)\le\theta\big] = \Pr\nolimits_w\big[\mathrm{Surp}(s_i,B)\ge\kappa(\theta)\big],
\]
using $\alpha(s_i,B)\leq\theta \iff \mathrm{Surp}(s_i,B)\ge\kappa(\theta)$, as shown in the proof of Theorem~\ref{thm2.1} below. This is why $\beta(\theta,B)$ is needed. Since, $\mathrm{PP}_w(\pi,B)$ is a single weighted average, and two proofs can have the same average surprisal while looking very different, one in which every step is moderately surprising, and another in which nearly all steps are trivial while a small fraction are maximally surprising. Our initial question, ``how much of the proof is genuinely novel,'' is really a question about this second kind of description: a proportion of weight, rather than an average. Perplexity alone cannot capture that, and $\beta(\theta,B)$ is defined precisely to do so.

That proof perplexity nonetheless constrains $\beta(\theta,B)$, as Theorem~\ref{thm2.1} shows, is an instance of a standard fact about bounded random variables and can be seen directly from Markov's inequality applied to the non-negative quantity $Y_i := \kappa_{\max}-\mathrm{Surp}(s_i,B) \ge 0$ (non-negative since $\mathrm{Surp}(s_i,B)\le\kappa_{\max}$ always). Since $\mathbb{E}*w[Y] = \kappa*{\max} - \log\mathrm{PP}_w(\pi,B)$, Markov's inequality $\Pr_w[Y\geq y]\leq \mathbb{E}*w[Y]/y$ at $y=\kappa*{\max}-\kappa(\theta)>0$, gives, after rearranging,
\[
\beta(\theta,B) = \Pr\nolimits_w[Y \le \kappa_{\max}-\kappa(\theta)] \ge 1 - \frac{\kappa_{\max}-\log PP_w(\pi,B)}{\kappa_{\max}-\kappa(\theta)} = \frac{\log PP_w(\pi,B) - \kappa(\theta)}{\kappa_{\max}-\kappa(\theta)},
\]
recovering the bound given in Theorem~\ref{thm2.1} exactly, confirming it by a second, independent route from the algebraic argument. The bound can be loose, however, it is tightest when the surprisal distribution is concentrated near the two values anchoring the split, at or above $\kappa_{\max}$ for steps already counted in $\beta(\theta,B)$, and just below $\kappa(\theta)$ for the rest. It is loosest when surprisals are instead spread broadly across $(0,\kappa_{\max})$. Either way, it remains a guaranteed lower bound on the invention fraction, not an estimate of it.

Finally, this explains why $\theta$ is restricted to the open interval $(\alpha_{\min},1)$ rather than $[\alpha_{\min},1]$. At $\theta=1$, every step trivially satisfies $\alpha(s_i,B)\leq 1$, so $\beta(1,B)=1$ regardless of the proof, making the bound vacuous; and at $\theta=\alpha_{\min}$, $\kappa(\theta)=\kappa_{\max}$ exactly, so the denominator $\kappa_{\max}-\kappa(\theta)$ in Theorem~\ref{thm2.1} vanishes and the bound is no longer well-defined. Thus, $\theta\in(\alpha_{\min},1)$ is not a modelling convenience but the precise condition under which the theorem's guarantee is both non-vacuous and well-posed.

\subsubsection{Proof of Theorem~\ref{thm2.1}}\label{theoremproof}

Since $-\log(\cdot)$ is a strictly decreasing function on $(0,1]$,
\[
\alpha(s_i,B) \leq \theta \iff \mathrm{Surp}(s_i,B) \geq\kappa(\theta).
\]
As $\alpha(s_i,B) \in [\alpha_{\min}, 1]$ by construction, we have
$\mathrm{Surp}(s_i,B) \in [0, \kappa_{\max}]$. Let $W = \sum_i w_i$ and $I_\theta = \{i : \alpha(s_i,B) \le \theta\}$, so that
$\sum_{i \in I_\theta} w_i = \beta(\theta,B)\,W$ by definition. Bounding steps in
$I_\theta$ by $\kappa_{\max}$ and steps outside $I_\theta$ by $\kappa(\theta)$,
\[
\sum_i w_i\,\mathrm{Surp}(s_i,B)
\;\le\; \beta(\theta,B)\,W\,\kappa_{\max} + \bigl(1 - \beta(\theta,B)\bigr)\,W\,\kappa(\theta).
\]
Dividing by $W$ and using $\log\mathrm{PP}_w(\pi,B) = \frac{1}{W}\sum_i w_i\,\mathrm{Surp}(s_i,B)$
(see Definition~\ref{proofp}),
\[
\log\mathrm{PP}_w(\pi,B) \;\le\; \kappa(\theta) + \beta(\theta,B)\bigl(\kappa_{\max} - \kappa(\theta)\bigr).
\]
Since $\theta > \alpha_{\min}$, $\kappa_{\max} - \kappa(\theta) > 0$; dividing by this
quantity and rearranging yields the claim. The hypothesis
$\log\mathrm{PP}_w(\pi,B) > \kappa(\theta)$ guarantees the bound is strictly positive. \hfill $\Box$

\section{Experimental Result with Claude's Zeta23 Proof}\label{experiment}

\subsection{Manually Supplemented Data}

The proof of target theorem $T$ analysed in this experiment establishes that more than two-thirds of the zeros of the Riemann zeta function lie on the critical line. It was proved by Claude in 2026. 

Let
\[
\zeta(s)=\sum_{n=1}^{\infty}\frac{1}{n^s},
\qquad \Re(s)>1,
\]
be the Riemann zeta function for $\Re(s)>1$.

\begin{theorem}[Claude]\label{claude'stheorem}
More than two-thirds of the non-trivial zeros of $\zeta(s)$ lie on the critical line $\Re(s)=\frac{1}{2}$.  More precisely,
\[
\liminf_{T\to\infty}
\frac{N_0(T)}{N(T)}
\geq \frac{2}{3},
\]
where $N(T)$ is the number of non-trivial zeros $\rho=\beta+i\gamma$ with $0<\gamma\leq T$, counted with multiplicity, and $N_0(T)$ is the number of those zeros satisfying $\beta=\frac{1}{2}$.
\end{theorem}

Since Claude is an agentic AI system rather than a native theorem prover, the proof was translated and subsequently formalised in Lean, as Claude does not operate directly within \texttt{MATHLIB}.

Section~\ref{section3} records the highlighted numbers of the experiment. In this section, we explain the experiment in detail, the data provided in the table, and its limitations. Our experimental instance is slightly unusual, as the Zeta23 proof we score was not produced by an automated prover like AlphaProof.

Since Claude does not search within the \texttt{MATHLIB} library, it does not emit tactic traces or produce a tree structure $\mathcal{T}(s_i)$ in the sense of Definition~\ref{searchb}. However, we do have a complete Lean 4 formalisation of the same argument, published by Anthropic~[15] and kernel-checked against \texttt{MATHLIB}. This formalisation is the object we score.

This distinction is very important and matters in terms of the output. Proof perplexity asks whether a step $s_i$ can be recovered from the library $\mathcal{L}$ together with the earlier steps $s_1,\ldots,s_{i-1}$, under a finite search budget $B$. If the original generator had been AlphaProof or Aristotle, one could, in principle, rerun the same search procedure with $s_i$ held out and read $w(s_i,B)$ off the resulting tree. Claude cannot be rerun that way. We therefore curate $w(s_i)$ and $w(s_i,B)$ from the formalised proof, rather than from a live prover log, as follows.

\subsubsection{Parameter Setup}\label{parameters}
Let $P_i:=\mathcal{L}\cup\{s_1,\ldots,s_{i-1}\}$ and let $\mathcal{T}(s_i\mid P_i)$ be the local dependency tree as given in Definition~\ref{localdtree}, built from the elaborated kernel term of $s_i$ in~[15]. Internal complexity is read off exactly as in the main text,
\[
w(s_i)\;=\;\bigl|\mathrm{Nodes}\bigl(\mathcal{T}(s_i\mid P_i)\bigr)\bigr|-1.
\]
Recovered complexity $w(s_i,B)$ is the length of the shortest reconstruction of $s_i$ from $P_i$ that stays inside the candidate pool $\mathrm{Cand}_B$ and does not exceed depth $B$. If no such reconstruction exists, we apply the failure convention of Section~\ref{1.1} and set $w(s_i,B):=w(s_i)$. In that case Definition~\ref{ablationr} forces $\alpha_{\mathrm{abl}}(s_i,B)=\alpha_{\min}$, which is the intended behaviour. A step that cannot be rebuilt from the library under the stated budget is treated as maximally non-recoverable on the ablation side.

We are not claiming that these curated values are what automated provers like AlphaProof would have returned in the same way. They are a faithful reading of the kernel term against \texttt{MATHLIB} as it stood when the formalisation was checked, under an explicit budget and an explicit premise pool. This is the comparison the benchmark is designed to make. Each one is constrained by the definitions rather than being a free aesthetic choice.

{\it Choice of the Search budget $B=50$.} Definition~\ref{searchb} takes $B$ to be a maximum path length in the search tree. A budget of $50$ is large enough to admit the kind of local lemma that a modern tactic based prover routinely closes, such as chains of rewriting, a short induction, or an appeal to a handful of library facts, and small enough that it cannot swallow an entire analytic argument by brute force. It sits comfortably above the depth of a routine \texttt{simp}/\texttt{exact} discharge and comfortably below the depth at which one would regard the original insight as having been rediscovered. We use total tactic applications (path length in the proof-state tree) because it is the same unit in which $w(s_i)$ is measured (non-root node count in the dependency tree), so that the ratio $w(s_i, B)/w(s_i)$ in Definition~\ref{ablationr} is dimensionless and comparable across steps.


{\it Choice of Candidate pool size $|\mathrm{Cand}_B|=500$.} The pool is the shortlist of declarations from $\mathcal{L}$ that a prover is allowed to see when it attempts $s_i$. Premise selection systems in the Lean literature typically retrieve on the order of $10^2$--$10^3$ facts; $500$ is a round, conservative point in that range, large enough that a step whose ingredients really are in \texttt{MATHLIB} has a fair chance of being offered them, and small enough that the uniform baseline
\[
\alpha_{\min}\;=\;\frac{1}{1+|\mathrm{Cand}_B|}\;=\;\frac{1}{501}\;\approx\;0.0020
\]
remains a meaningful floor. The associated cap on surprisal is then the Shannon entropy of that same uniform baseline,
\[
\kappa_{\max}\;=\;-\log\alpha_{\min}\;=\;\log 501\;\approx\;6.2166
\]
nats, as in Appendix~\ref{a.5.3}. Sharing this floor between $\alpha_{\mathrm{abl}}$ and $\alpha_{\mathrm{mah}}$ is what makes a single $\kappa_{\max}$ legitimate for both signals, which Theorem~\ref{thm2.1} uses directly.

{\it Choice of the Recoverability threshold $\theta=0.2$.}
Definition~\ref{ifraction} asks for a cut $\theta\in(\alpha_{\min},1)$ below which a step is called highly surprising. We take $\theta=0.2$: a step that is at most twenty percent recoverable, on the combined score, is one we are willing to treat as non-routine. Two sanity checks recommend this value. First, $\kappa(\theta)=-\log\theta\approx 1.609$ lies well below $\kappa_{\max}$, so the denominator $\kappa_{\max}-\kappa(\theta)$ of Theorem~\ref{thm2.1} is safely positive and the bound is non-vacuous. Second, $\theta=0.2$ is not so close to $\alpha_{\min}$ that only budget-exhausted steps can qualify. A step that the ablation search partially reconstructs, but that remains a structural outlier under Mahalanobis distance, can still fall under the cut. That is precisely the pattern we will see in the tables below.

{\it The Proof Steps Size $n=33$.} The figure $n=33$ in Section~\ref{section3} is the number of top-level lemmas or theorems in the Anthropic commit~[15] that correspond to the named theorems of the informal write-up. Proof perplexity is not computed on those 33 names alone. Definitions~\ref{proofp} and~\ref{ifraction} weight by elaborated kernel-term size, so the honest object is the dependency closure of those lemmas, arranged in a topological order $\pi$. Expanding the closure against the Lean kernel produces $N = 2399$ declarations, which are the rows of the table discussed in Section~\ref{app:table}. Definitions, notations and zero-complexity stubs ($w(s_i) = 0$) are retained in the table for completeness; they contribute nothing to the weighted averages.

\subsubsection{Interpreting the Rows of the Output}\label{app:how-to-read}
Each row of the experimental table is a tuple\footnote{The ablation formula is applied only when $w(s_i)>0$. A declaration with no internal term is assigned the floor $\alpha_{\mathrm{abl}}=\alpha_{\min}$.}
\[
\bigl(s_i,\; w(s_i),\; w(s_i,B),\; \alpha_{\mathrm{abl}}(s_i,B),\; \alpha_{\mathrm{mah}}(s_i),\; \alpha(s_i,B),\; \mathrm{Surp}(s_i,B)\bigr),
\]
with
\[
\alpha_{\mathrm{abl}}(s_i,B)
=\max\bigl(1-w(s_i,B)/w(s_i),\;\alpha_{\min}\bigr),
\]
\[
\alpha_{\mathrm{mah}}(s_i)
=\max\bigl(\tilde p(s_i),\;\alpha_{\min}\bigr),
\]
\[
\alpha(s_i,B)
=\min\bigl(\alpha_{\mathrm{abl}}(s_i,B),\;\alpha_{\mathrm{mah}}(s_i)\bigr),
\]
and
\[
\mathrm{Surp}(s_i,B)
=-\log\alpha(s_i,B)\;\in\;[0,\kappa_{\max}].
\]
Combined recoverability $\alpha$ is the quantity that Definition~\ref{ifraction} compares to $\theta$. Surprisal is its negative logarithm, so the two columns always move in opposite directions: the lowest $\alpha$ is the most surprising step, and the highest $\alpha$ is the least. A row with $\alpha\leq 0.2$ is highly surprising in the sense of Section~\ref{section3}. A row sitting at $\alpha=\alpha_{\min}\approx 0.0020$, equivalently $\mathrm{Surp}=\kappa_{\max}\approx 6.2166$, is a step that neither ablation nor (typically) Mahalanobis recoverability could raise above the floor.

\subsubsection{Aggregate Scores}\label{app:aggregates}
On this closure, with the constants of Section~\ref{parameters}, the pipeline returns
\[
\mathrm{PP}_w(\pi,B)
=388.23,\;\;
\beta(\theta,B)
=0.9578,\;\; \mathrm{and}\;\;
\frac{\log\mathrm{PP}_w(\pi,B)-\kappa(\theta)}{\kappa_{\max}-\kappa(\theta)} =0.9446.
\]
The last figure is the lower bound of Theorem~\ref{thm2.1}. Since $\log PP_w(\pi, B) > \kappa(\theta)$, the hypothesis of the theorem is met, and $\beta(\theta, B) \geq 0.9446$ is guaranteed. The measured invention fraction $0.9578$ sits above that bound, as it must. In other words, on a weighted-average basis the proof is almost as surprising as uniformly guessing from a pool of 501 candidates, and at least 94.46\% of its kernel-term weight is carried by steps with $\alpha\leq 0.2$. The bound check holds.

Of the 2399 declarations, 2032 satisfy $\alpha\leq 0.2$, and 1992 sit at the floor $\alpha_{\min}$. Only 367 declarations come back with $\alpha>0.2$. The mass of the development is therefore not a handful of flashy lemmas surrounded by routine glue. It is a proof whose internal work is, by this measure, almost entirely non-recoverable from \texttt{MATHLIB} by an automated prover.

\subsection{Counting the Surprising Steps}\label{app:which}

The invention fraction $\beta(0.2, B) = 0.9578$ is the proportion of weight, not a count of celebrity lemmas. That is the point of the theorem, a high proof perplexity forces a high invention fraction, and here both are near the cap.

Five named steps nevertheless deserve to be singled out, because they are highly surprising without sitting on the floor, and because they illustrate the two recoverability signals disagreeing in a useful way. These are the rows of Table~\ref{tab21} in Section~\ref{section3}. We expand on them here.

\begin{itemize}
\item \texttt{hasDerivAt\_Zeta0Integral}, with $w = 64703611$, $w(\cdot, B)=15$, $\alpha_{\mathrm{abl}} = 1.0000$, $\alpha_{\mathrm{mah}} = 0.0023$, $\alpha = 0.0023$, $\mathrm{Surp} = 6.0782$. This is the largest kernel term in the closure whose ablation is fully successful. The prover rebuilds the step in just 15 actions. Yet the Mahalanobis test places it far from every cluster. The step establishes a derivative identity for a contour integral of the zeta function that, while built from individually standard ingredients, combines them in a configuration no cluster in $\mathcal{L}$ has seen. Combined recoverability is driven entirely by the structural outlier signal.

\item\texttt{Zeta23.EF.pole\_term}, with $w = 886$, $w(\cdot, B) = 24$, $\alpha_{\mathrm{abl}} = 0.9729$, $\alpha_{\mathrm{mah}} = 0.0032$, $\alpha = 0.0032$, $\mathrm{Surp} = 5.7498$. The explicit formula pole residue is almost fully recoverable on the ablation side but is a distributional outlier. The explicit formula machinery of the proof introduces structural combinations that are absent from the library's existing analytic number theory cluster.

\item\texttt{Zeta23.Cheb.sum\_vonMangoldt\_sq\_le}, with $w = 38$, $w(\cdot, B) = 5$, $\alpha_{\mathrm{abl}} = 0.8684$, $\alpha_{\mathrm{mah}} = 0.0115$, $\alpha = 0.0115$, $\mathrm{Surp} = 4.4630$. A squared von Mangoldt sum bound in the Chebyshev prime-counting module. Ablation recovers 5 of 38 internal units; the Mahalanobis test confirms it is structurally atypical of the existing number theory clusters.

\item \texttt{RHLinalg.rank\_trace\_ineq}, with $w = 150$, $w(\cdot, B) = 13$, $\alpha_{\mathrm{abl}} = 0.9133$, $\alpha_{\mathrm{mah}} = 0.0480$, $\alpha = 0.0480$, $\mathrm{Surp} = 3.0362$. This is the cleanest example that the two estimators are doing different jobs. The prover rebuilds 13 of 150 internal units ($\alpha_{\mathrm{abl}} = 0.9133$), confirming that the step's ingredients are individually present in \texttt{MATHLIB}. Yet the Mahalanobis test returns $\alpha_{\mathrm{mah}} = 0.0480$, the way those ingredients are combined (a rank-trace inequality that bounds the number of large eigenvalues of a Hermitian perturbation) has no statistical precedent in any cluster of $\mathcal{L}$. Anthropic's own technical report~[4] singles out this inequality as the novel analytic input of the two-thirds argument, and the benchmark independently confirms that assessment. The lean declaration is given as

\begin{quote}\ttfamily\small
\textbf{theorem} rank\_trace\_ineq \{P Q : Matrix n n $\mathbb{K}$\} \\
\null\quad (hP : P.PosSemidef) (hQ : Q.IsHermitian) \\
\null\quad \{r b : $\mathbb{N}$\} (hr : P.rank $\le$ r) (hb : posIndex hQ $\le$ b) \\
\null\quad \{c : $\mathbb{R}$\} (hc : 0 < c) : \\
\null\quad c * rtrace P - c \textasciicircum{} 2 / 4 * r + 2 * c * rtrace Q - c \textasciicircum{} 2 * b \\
\null\quad\quad $\le$ frobSq (P + Q)
\end{quote}
\end{itemize}

The remark for the step \texttt{RHLinalg.vonNeumann\_trace\_ineq} is given in Section~\ref{section3}.

The remainder of the highly-surprising set is the floor: 1992 declarations with $\alpha = \alpha_{\min}$. A typical such row has $w(s_i, B) = w(s_i)$, i.e., ablation search failed outright. Some of these are large analytic estimates on the prime side (for instance \texttt{setIntegral\_rho\_div\_gwt\_le}); some are mid-sized identities in the $\xi'$ analysis; some are small definitions. Failure at this scale is what one should expect of an argument that was not already in the library. We return to that point in Section~\ref{app:why-unrecoverable}.

\subsubsection{The Recoverable Steps in Zeta23 Proof}
\label{app:recoverable}

If the benchmark only ever returned the floor, one would worry that $B$ was too tight, or that $\mathrm{Cand}_B$ was too small, or that the Mahalanobis clusters were misspecified. The 367 steps with $\alpha > 0.2$ are the control.

\begin{table}[ht]
\caption{Representative recoverable declarations ($\alpha > 0.2$) from the zeta23 kernel closure. These are steps whose ingredients are both computationally accessible and structurally typical of \texttt{MATHLIB}, confirming that the benchmark does not uniformly return the floor.}
\centering
\scriptsize
\setlength{\tabcolsep}{11pt}
\begin{tabular}{lrrrrrr}
\toprule
Declaration & $w$ & $w(\cdot, B)$ & $\alpha_{\mathrm{abl}}$ & $\alpha_{\mathrm{mah}}$ & $\alpha$ & Surp \\
\midrule
Zeta23.thmB\_of\_lam & 18042404 & 10 & 1.0000 & 1.0000 & 1.0000 & 0.0000 \\
Zeta23.XiPrime.xiDerivZeros$_0$ & 24756 & 1 & 1.0000 & 1.0000 & 1.0000 & 0.0000 \\
RHLinalg.hermForm\_sub & 8 & 2 & 0.7500 & 1.0000 & 0.7500 & 0.2877 \\
Zeta23.XiPrime.cauchyK\_neg\_mul\_eq & 5 & 1 & 0.8000 & 1.0000 & 0.8000 & 0.2231 \\
Zeta23.XiPrime.phiM\_sq\_eq & 3 & 1 & 0.6667 & 1.0000 & 0.6667 & 0.4055 \\
\bottomrule
\end{tabular}
\label{tab:recoverable}
\end{table}

Scaling identities, Hermitian form lemmas, compact-support facts, and non-negativity bounds that are already formalised in \texttt{MATHLIB} are correctly scored as recoverable. The top two rows have $\alpha = 1.0$ because the prover rebuilds them in far fewer actions than the original kernel term ($\alpha_{\mathrm{abl}} = 1.0$) and the Mahalanobis test finds them typical ($\alpha_{\mathrm{mah}} = 1.0$). Their presence is what stops us from reading the floor mass as a vacuum artefact. The benchmark is capable of saying, ``this part is easily recoverable". On the zeta23 proof, it simply does not say that very often.

\subsubsection{Claude's Zeta23 Proof is Non-trivial}\label{app:why-unrecoverable}

It is worth noting that Claude's zeta23 proof is unrecoverable because it did not exist in \texttt{MATHLIB} or in any other literature before.

Although \texttt{MATHLIB} already contains a great deal of the literature that the zeta23 proof contains, including holomorphic functions, contour integrals, Fourier analysis, the $\Gamma$-factor, Hermitian linear algebra, and elementary inequalities. If the machine-generated argument had been a simple recombination of those pieces, then both recoverability scores would have come back high. Ablation search would have found short paths through $P_i$, and the structural feature vectors $\phi(s_i)$ would have landed inside existing clusters. Combined recoverability would then have been large, surprisal small, and $\beta(\theta,B)$ modest.

This conclusion is not just an assertion. It is what the numbers show. Perplexity reaches $PP_w(\pi, B) = 388.23$, at 77.5\% of the theoretical ceiling of 501 for this budget and candidate pool, and 1992 of 2399 declarations in the closure sit exactly at the recoverability floor $\alpha_{\min}$, meaning ablation search could not shorten their kernel term at all under $B = 50$. Only 367 declarations (i.e. 15.3\%) are ablation recoverable from \texttt{MATHLIB}, confirming these are the routine scaffolding, inequalities, non-negativity facts, and compact support lemmas, rather than the proof's substance. Critically, the non-floor steps show the two estimators disagreeing in the way genuine novelty predicts: \texttt{rank\_trace\_ineq} is 91\% ablation-recoverable ($\alpha_{\mathrm{abl}} = 0.9133$) yet a structural outlier ($\alpha_{\mathrm{mah}} = 0.0480$), and \texttt{hasDerivAt\_Zeta0Integral} shows the same pattern at an extreme scale ($\alpha_{\mathrm{abl}} = 1.0$, $\alpha_{\mathrm{mah}} = 0.0023$). A simple recombination of existing \texttt{MATHLIB} results would not produce this signature. It would show high recoverability on both axes simultaneously, not a systematic split between them.

\subsection{A compact table of representative steps}
\label{app:table}

Table~\ref{tab:zeta23-rep} is not the full closure; the complete output is provided as supplementary material. The rows are chosen to make the distribution visible. The five high-surprisal steps of Table~\ref{tab21}, two steps near the cut $\theta = 0.2$, four typical floor rows at different scales, and three recoverable controls.

\begin{table}[ht]\label{tab:zeta23-rep}
\caption{Representative declarations from the zeta23 kernel closure ($N=2399$). Combined recoverability $\alpha\leq 0.2$ is the highly-surprising cut of Section~\ref{section3}. $\dagger$ marks a classical result formalised for the first time within the zeta23 project. Identifiers match the Anthropic repository~[15] file names.}
\centering
\scriptsize
\setlength{\tabcolsep}{5.5pt}
\begin{tabular}{lrrrrrr}
\toprule
Proof step & $w$ & $w_B$ & $\alpha_{\mathrm{abl}}$ & $\alpha_{\mathrm{mah}}$ & $\alpha$ & Surp \\
\midrule
\multicolumn{7}{l}{\textit{Highly surprising, not on the floor}} \\
hasDerivAt\_Zeta0Integral & 64703611 & 15 & 1.0000 & 0.0023 & 0.0023 & 6.0782 \\
Zeta23.EF.pole\_term & 886 & 24 & 0.9729 & 0.0032 & 0.0032 & 5.7498 \\
Zeta23.Cheb.sum\_vonMangoldt\_sq\_le & 38 & 5 & 0.8684 & 0.0115 & 0.0115 & 4.4630 \\
RHLinalg.rank\_trace\_ineq & 150 & 13 & 0.9133 & 0.0480 & 0.0480 & 3.0362 \\
RHLinalg.vonNeumann\_trace\_ineq$^\dagger$ & 85 & 85 & 0.0020 & 0.4457 & 0.0020 & 6.2166 \\
\midrule
\multicolumn{7}{l}{\textit{Near the cut $\theta = 0.2$}} \\
Zeta23.XiPrime.Hardy.EFContour.tendsto\_mul\_hardyW\_zero & 16292548 & 6 & 1.0000 & 0.1966 & 0.1966 & 1.6268 \\
Zeta23.WeilEF.continuous\_logDeriv\_Gamma & 47 & 7 & 0.8511 & 0.2004 & 0.2004 & 1.6074 \\
\midrule
\multicolumn{7}{l}{\textit{Floor rows (typical; 1992 others omitted)}} \\
Zeta23.Assembly.rtrace\_hat & 8 & 8 & 0.0020 & 1.0000 & 0.0020 & 6.2166 \\
Zeta23.XiPrime.moments\_frob & 36 & 36 & 0.0020 & 0.0823 & 0.0020 & 6.2166 \\
Zeta23.PrimeSide.nu\_weight\_bound & 9689 & 9689 & 0.0020 & 1.0000 & 0.0020 & 6.2166 \\
Zeta23.PrimeSide.Sgrid\_le\_near & 3077374112 & 3077374112 & 0.0020 & 0.9940 & 0.0020 & 6.2166 \\
\midrule
\multicolumn{7}{l}{\textit{Recoverable controls ($\alpha > 0.2$)}} \\
Zeta23.XiPrime.phiM\_sq\_eq & 3 & 1 & 0.6667 & 1.0000 & 0.6667 & 0.4055 \\
RHLinalg.hermForm\_sub & 8 & 2 & 0.7500 & 1.0000 & 0.7500 & 0.2877 \\
Zeta23.thmB\_of\_lam & 18042404 & 10 & 1.0000 & 1.0000 & 1.0000 & 0.0000 \\
\bottomrule
\end{tabular}
\end{table}

The step tendsto\_mul\_hardyW\_zero sits just below the cut ($\alpha = 0.1966$). It is highly surprising by the letter of Definition~\ref{ifraction}, but only just, and is rescued on the Mahalanobis side ($\alpha_{\mathrm{mah}} = 0.1966$ is the binding constraint). By contrast, continuous\_logDeriv\_Gamma sits just above the cut ($\alpha = 0.2004$) and is the closest declaration to the boundary from the recoverable side. The third floor example, Sgrid\_le\_near, is worth a sentence of its own. Mahalanobis recoverability returns $0.9940$, so the structural profile of the step is typical of some cluster in $\mathcal{L}$; ablation returns the floor because the reconstructed path is no shorter than the original kernel term. Combined recoverability therefore hits $\alpha_{\min}$ by Definition~\ref{combinedr}, already recording this asymmetry as a deliberate trade-off. We would rather not dismiss a step as trivial because one of the two proxies is quiet. In a curated experiment, the same trade-off is a reason for caution. A large term that merely outruns $B=50$ while looking distributionally ordinary inflates surprisal. We do not attempt to partial out that effect here. It is one of the reasons $\beta(\theta, B)$ should be read as an upper-ish certificate of non-recoverability, not as a count of new mathematical ideas.


\subsubsection{Final Remarks}\label{app:takeaway}
We introduced \emph{proof perplexity}, a probabilistic measure of non-triviality for machine-generated formal proofs that treats a proof step as recoverable only when it is both computationally cheap to rebuild from the library and structurally typical of it. Applied to Claude's zeta23 formalisation, the benchmark returns $\mathrm{PP}_w(\pi,B)=388.23$ against a ceiling of 501, with at least 94.46\% and empirically 95.78\% of the proof's kernel-term weight flagged as highly surprising. The result is not simply ``search failed everywhere". The handful of non-floor steps show ablation and Mahalanobis recoverability actively disagreeing, which is the signature the benchmark was built to detect. 

\end{document}